\documentclass{amsart}
\usepackage{graphicx}
\usepackage[T1]{fontenc}
\usepackage{amsmath}
\usepackage{amssymb}
\usepackage{amsthm}
\usepackage{epsfig}
\usepackage{color}
\usepackage{hyperref}
\usepackage{enumitem}
\usepackage{listings}
\usepackage{lscape}
\definecolor{gray}{rgb}{0.5,0.5,0.5}

\newlist{steps}{enumerate}{1}
\setlist[steps, 1]{label = Step \arabic*:}

\input xy
\xyoption{all}
\usepackage[all]{xy}

\newtheorem{theorem}{Theorem}[section]

\theoremstyle{definition}
\newtheorem{definition}[theorem]{Definition}
\newtheorem{example}[theorem]{Example}

\theoremstyle{remark}
\newtheorem{remark}[theorem]{Remark}

\def\MR#1{\href{http://www.ams.org/mathscinet-getitem?mr=#1}{MR}}
\DeclareMathOperator{\fix}{Fix}

\begin{document}

\title[Lights Out on a Random Graph]{Lights Out on a Random Graph}

\author{Bradley Forrest}
\address{101 Vera King Farris Drive; Galloway, NJ 08205}
\email{bradley.forrest@stockton.edu}

\author{Nicole Manno}
\email{nicoleannemanno@gmail.com}



\begin{abstract}
    We consider the generalized game \emph{Lights Out} played on a graph and investigate the following question: for a given positive integer $n$, what is the probability that a graph chosen uniformly at random from the set of graphs with $n$ vertices yields a universally solvable game of \emph{Lights Out}?  When $n \leq 11$, we compute this probability exactly by determining if the game is universally solvable for each graph with $n$ vertices. We approximate this probability for each positive integer $n$ with $n \leq 100$ by applying a Monte Carlo method using 1,000,000 trials.  We also perform the analogous computations for connected graphs.
\end{abstract}

\maketitle

\section{Introduction}

\emph{Lights Out} is an electronic one-player game played on a $5 \times 5$ square grid of lights.  Each one of the 25 square lights, when pressed, toggles itself and each of the lights with which it shares an edge. The goal in \emph{Lights Out} is to produce a configuration in which all of the lights are turned off.

In this work, we consider \emph{Lights Out} played on an undirected graph $\Gamma = (V,E)$, where there is one light for each element of $V$ and pressing $v \in V$ toggles $v$ and each vertex adjacent to $v$.  We will call $C \subseteq V$ the \emph{initial configuration} of a game of \emph{Lights Out} on $\Gamma$ where $C$ is the set of vertices that are toggled on when the game begins.  Winning \emph{Lights Out} consists of finding a subset of $V$ so that pressing the lights in that subset toggles exactly $C$; indeed it is only necessary to find a subset of $V$ and not a sequence of that subset since the effect of pressing lights is independent of the order in which they are pressed.  If such a subset of $V$ exists, then we say that the game is \emph{solvable}.  Much is known about the solvability of particular initial configurations. For any graph $\Gamma$, if $C = V$ then the game is solvable, and if $V$ is an $n \times n$ square grid with $n$ odd and $C$ is only the center square, then the game is solvable \cite{CHKR} \cite{EES} \cite{HY} \cite{Sutner1}.  We say that \emph{Lights Out} on $\Gamma$ is \emph{universally solvable} if it is solvable for every initial configuration $C$. For example, the $3 \times 3$ square grid is universally solvable while the $5 \times 5$ grid is not \cite{Sutner1}.

In \cite{AS1} and \cite{AS2}, Amin and Slater made substantial progress determining which graphs are universally solvable.  In particular, they not only give equivalent conditions to determine when a graph is universally solvable but also classify the paths, spiders, and caterpillars that are universally solvable and provide a method to generate all universally solvable trees.  The existence of so many classes of graphs that are universally solvable raises the primary question considered in this work: if a simple graph is chosen uniformly at random from the set of graphs with $n$ vertices, what is the probability that it will yield a universally solvable game of \emph{Lights Out}?  Since games of \emph{Lights Out} on disconnected graphs can be thought of as independent games on connected graphs, we also consider the analogous question for connected graphs. To approach this question, we implemented the algorithm for choosing a graph uniformly at random described by Dixon and Wilf in \cite{Dixon}.  Our program was written in Java and our code is available at \url{http://github.com/nicolemanno/Lights-Out}.

This paper is organized as follows.  In Section \ref{sec:lin}, we discuss the connection between \emph{Lights Out} and linear algebra and give, for each $n \leq 11$, the number of graphs with $n$ vertices that are universally solvable. We present the algorithm that we used to select a graph uniformly at random in Section \ref{sec:random} and give the results of our Monte Carlo experiments.  In Section \ref{sec:future}, we discuss possible extensions of our work and interesting open problems.

\section{Lights Out and Linear Algebra}
\label{sec:lin} 

\textit{Lights Out} can be studied through linear algebra by ordering the vertices of $\Gamma$.  Once an order is chosen, subsets of the vertices of $\Gamma$ can be represented by column vectors.  Let $\Gamma$ have $n$ vertices and for simplicity we number the vertices from 1 to $n$. The subset $C$ is given by the $n \times 1$ column vector 
\[ \overrightarrow{b} = \left[
\begin{array}{cccc}
b_{1} & b_{2} & \cdots & b_{n}
\end{array} \right]^T, \]
where $b_{i} = 1$ if vertex $i$ is in $C$ and is 0 otherwise.

The primary tool that we will use to study \emph{Lights Out} is the \emph{neighborhood matrix} of the graph $\Gamma$:

\begin{definition}
Let $\Gamma = (V, E)$ be a graph with $n$ vertices labeled $1$ to $n$. The \emph{neighborhood matrix} $\mathcal{A} = [ a_{i,j} ]$ of $\Gamma$ is the $n \times n$ matrix in which $a_{i,j} = 1$ if $i=j$ or if the $i$-th and $j$-th vertices are adjacent and $ a_{i,j}=$ 0 otherwise. 
\end{definition}

\begin{remark}
The matrix $\mathcal{A}$ is symmetric and the $i$-th row of $\mathcal{A}$ has a 1 in the $j$-th column if and only if pressing vertex $i$ toggles vertex $j$. This matrix can also be expressed as the sum of the adjacency matrix of $\Gamma$ and the $n \times n$ identity matrix.
\end{remark}

When the column vector $\overrightarrow{b}$ represents the initial configuration of \emph{Lights Out} on a graph $\Gamma$ whose vertices are labeled $1$ to $n$, solving the game reduces to finding the vector $\overrightarrow{x}$ so that $\mathcal{A}\overrightarrow{x} = \overrightarrow{b},$ where all of the coefficients are in the field $\mathbb{Z}_2$.
The initial configuration represented by $\overrightarrow{b}$ is solvable if and only if it belongs to the column space of $\mathcal{A}$ \cite{Anderson}. This leads to the following fundamental observation about \emph{Lights Out} on $\Gamma$: \cite{GKT} \cite{GKTZ} \cite{Sutner1} \cite{Sutner2}

\begin{theorem} \label{thm:main}
Lights Out on a graph $\Gamma$ is universally solvable if and only if the neighborhood matrix of $\Gamma$ is invertible.
\end{theorem}

When $\mathcal{A}$ is invertible, winning the game requires pressing the vertices represented by $\overrightarrow{x} = \mathcal{A}^{-1}\overrightarrow{b}$.  Our program checks for invertibiliy by performing row reduction. Note that the order chosen for the vertices of $\Gamma$ does not affect the invertibility of $\mathcal{A}$ since two different orderings give neighborhood matrices that are conjugate by permutation matrices.

\begin{example}
Figure \ref{fig:graphs} shows labeled representatives of the six distinct unlabeled connected graphs with four vertices.  What follows are the neighborhood matrices corresponding to these labeled graphs.

\begin{figure}[!htb]
\begin{center}
\begin{picture}(300,150)
\put(24,12){\includegraphics[scale = .7]{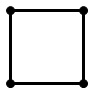}} 
\put(124,12){\includegraphics[scale = .7]{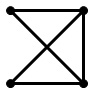}}
\put(224,12){\includegraphics[scale = .7]{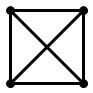}}
\put(24,87){\includegraphics[scale = .7]{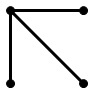}}
\put(124,87){\includegraphics[scale = .7]{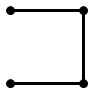}}
\put(224,87){\includegraphics[scale = .7]{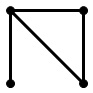}}

\put(19,55){1}
\put(19,130){1}
\put(119,55){1}
\put(119,130){1}
\put(219,55){1}
\put(219,130){1}

\put(73,55){2}
\put(73,130){2}
\put(173,55){2}
\put(173,130){2}
\put(273,55){2}
\put(273,130){2}

\put(73,10){3}
\put(73,85){3}
\put(173,10){3}
\put(173,85){3}
\put(273,10){3}
\put(273,85){3}

\put(19,10){4}
\put(19,85){4}
\put(119,10){4}
\put(119,85){4}
\put(219,10){4}
\put(219,85){4}

\end{picture}

\caption{\label{fig:graphs} Labeled representatives of the six unlabeled connected graphs with 4 vertices.}
\end{center}

\end{figure}
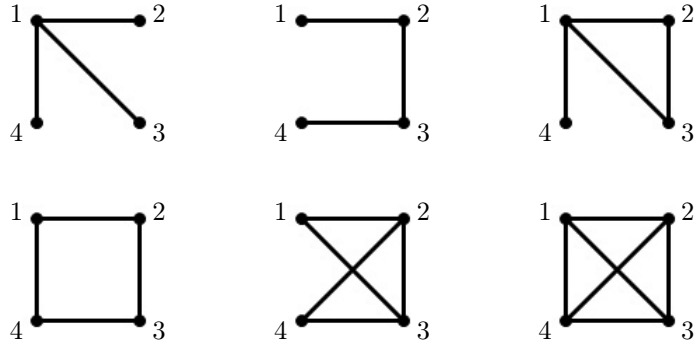

\[
\left[
\begin{array}{cccc}
1 & 1 & 1 & 1\\
1 & 1 & 0 & 0\\
1 & 0 & 1 & 0\\
1 & 0 & 0 & 1
\end{array}
\right]\;\;\;\;
\left[
\begin{array}{cccc}
1 & 1 & 0 & 0\\
1 & 1 & 1 & 0\\
0 & 1 & 1 & 1\\
0 & 0 & 1 & 1
\end{array}
\right]
\;\;\;\;
\left[
\begin{array}{cccc}
1 & 1 & 1 & 1\\
1 & 1 & 1 & 0\\
1 & 1 & 1 & 0\\
1 & 0 & 0 & 1
\end{array}
\right]
\]\\
\[
\left[
\begin{array}{cccc}
1 & 1 & 0 & 1\\
1 & 1 & 1 & 0\\
0 & 1 & 1 & 1\\
1 & 0 & 1 & 1
\end{array}
\right]
\;\;\;\;
\left[
\begin{array}{cccc}
1 & 1 & 1 & 0\\
1 & 1 & 1 & 1\\
1 & 1 & 1 & 1\\
0 & 1 & 1 & 1
\end{array}
\right]
\;\;\;\;
\left[
\begin{array}{cccc}
1 & 1 & 1 & 1\\
1 & 1 & 1 & 1\\
1 & 1 & 1 & 1\\
1 & 1 & 1 & 1
\end{array}
\right]
\] 

\bigskip

Row reducing each neighborhood matrix demonstrates that the only connected universally solvable graphs with 4 vertices are the four cycle and the length three path, represented in Figure \ref{fig:graphs} by the bottom left and top middle graphs, respectively.
\end{example}

Using the archives of graphs and connected graphs available from \cite{McKay1} and \cite{McKay2}, we applied Theorem \ref{thm:main} to determine which of these graphs correspond to universally solvable games of \emph{Lights Out}.  For 1  to 10 vertices, the program took less than 5 minutes to run, while the run time to compute results for graphs with 11 vertices was approximately 10 hours; these computations were executed on a 2015 MacBook Pro wtih 16 GB of RAM. The results for all graphs follow in Table \ref{tab:small} while the results for connected graphs are in Table \ref{tab:smallconn}.\\

\begin{table}[!htb]
\begin{center}
\caption{\label{tab:small} For each $n$ with $1 \leq n \leq 11$, the probability that a graph chosen uniformly at random from the set of graphs with $n$ vertices is universally solvable.}
\bigskip
\begin{tabular}{ccccccc}
\hline\\[-5 pt]
Number & & Number & & Number of & & Probability \\
of & & of & & Universally & & Universally\\
Vertices & & Graphs & & Solvable Graphs & & Solvable\\[5 pt]
\hline\\[-5 pt]
1 & & 1 & & 1  & &  1\\
2 & & 2 & & 1 & & 0.5\\
3 & & 4 & & 2  & & 0.5\\
4 & & 11 & & 4 & & 0.363636\\
5 & & 34 & & 13 & & 0.382353\\
6 & & 156 & & 47 & & 0.301282\\
7 & & 1044 & & 339 & & 0.324713\\
8 & & 12346 & & 4043 & & 0.327474\\
9 & & 274668 & & 98375 & & 0.358160\\
10 & & 12005168 & & 4553432 & & 0.379289\\
11 & & 1018997864 & & 403286335 & & 0.395768\\[5 pt]
\hline
\end{tabular}
\end{center}
\end{table}

\bigskip
\bigskip

\begin{table}[!htb]
\begin{center}
\caption{\label{tab:smallconn} For each $n$ with $1 \leq n \leq 11$, the probability that a graph chosen uniformly at random from the set of connected graphs with $n$ vertices is universally solvable.}
\bigskip
\begin{tabular}{ccccccc}
\hline\\[-5 pt]
Number & & Number of & & Universally & & Probability\\
of & & Connected & & Solvable & & Universally\\
Vertices & & Graphs & & Connected Graphs & & Solvable\\[5 pt]
\hline\\[-5 pt]
1 & & 1 & & 1 & &  1\\
2 & & 1 & &  0 & &  0\\
3 & & 2 & &  1 & &  0.5\\
4 & & 6 & &  2 & &  0.333333\\
5 & & 21 & &  9 & &  0.428571\\
6 & & 112 & &  33 & &  0.294643\\
7 & & 853 & &  290 & &  0.339977\\
8 & & 11117 & &  3692 & &  0.332104\\
9 & & 261080 & &  94280 & &  0.361115\\
10 & & 11716571 & &  4454654 & &  0.380201\\
11 & & 1006700565 & &  398728322 & &  0.396074\\[5 pt]
\hline
\end{tabular}
\end{center}
\end{table}

\section{Algorithm for Large Numbers of Vertices}
\label{sec:random}
When we consider graphs with more than 11 vertices, the methods of Section \ref{sec:lin} are ineffective; for instance there are more than $10^{11}$ unlabeled graphs with 12 vertices. For this reason we applied a Monte Carlo method.  For each $n$ with $1 \leq n \leq 100$, we ran 1,000,000 trials of the experiment of selecting a graph uniformly at random from the set of graphs with $n$ vertices and determined if each chosen graph produced a universally solvable game of \emph{Lights Out}.  To choose a graph uniformly at random, we followed the algorithm described by Dixon and Wilf in \cite{Dixon}.  We will review their methods below.

\subsection{Selection of a Graph Uniformly At Random}

 Consider the action of the symmetric group $S_n$ on the set of all labeled graphs with $n$ vertices where the action is given by permuting the labels.  The unlabeled graphs with $n$ vertices are in one-to-one correspondence with orbits of this action.  For $g \in S_n$, let $\fix(g)$ be the set of fixed points of $g$ under the action and consider the set $X = \{(g,\Gamma) \mid g \in S_n, \Gamma \in \fix(g) \}$. Dixon and Wilf observe that each orbit is represented the same number of times within $X$.  By choosing an element from each conjugacy class of $S_n$, each orbit is also represented the same number of times within the subset of $X$ corresponding to the chosen elements of $S_n$.  So, we can choose a graph uniformly at random by
 \begin{itemize}
     \item first choosing a conjugacy class of $S_n$ weighted by the product of the cardinality of its fixed point set and the number of elements in the conjugacy class,
     \item and then choosing a graph uniformly at random from the fixed point set of a representative of the chosen conjugacy class.
 \end{itemize}
 
 Note that the conjugacy classes of $S_n$ are in one-to-one correspondence with partitions of $n$.  More specifically, the partition $[k_1,k_2, \ldots, k_n]$ with $k_i$ parts of size $i$ corresponds to the conjugacy class containing permutations with $k_i$ cycles of length $i$ for each $i$ with $1 \leq i \leq n$.  So, we label conjugacy classes by their corresponding partitions.  Dixon and Wilf give the conjugacy class $[k_{1},\ldots ,k_{n}]$ weight $$w([k_{1},...,k_{n}]) = \frac{n!2^{c(g)}}{\displaystyle \prod_{i=1}^n(i^{k_{i}}k_{i}!)},$$ where $$c(g) = \frac{1}{2}\left(\displaystyle\sum_{i=1}^{n}l(i)^{2}\phi(i) - l(1)+l(2) \right), \qquad \qquad l(i) = \displaystyle\sum_{j=1, \, i|j}^n k_{j},$$ and $\phi$ is the Euler phi function.  Following their algorithm, we choose a partition $\pi$ of $n$ so that the probability of choosing $[k_{1},\ldots,k_{n}]$ is equal to $w([k_{1},\ldots,k_{n}])/n!g_{n}$, where $g_{n}$ is the number of unlabeled graphs with $n$ vertices. Our methods require knowledge of the value of $g_n$ and this value has been computed up to $n = 140$ by Briggs \cite{Briggs} \cite{oeis} using methods of Oberschelp \cite{Ober}.  While it is feasible to extend our results to $n = 140$, as $n$ increases our approximations stabilize to the narrow range from 0.4185 to 0.4210 which limits the utility of extending the computation.

We follow the suggestion of Dixon and Wilf \cite{Dixon} and apply results of Oberschelp \cite{Ober} to efficiently choose a partition. Within the set of partitions of $n$, the weight is concentrated in the partitions with the most parts of size 1. We choose a random number $\xi$, such that $0 \leq \xi < 1$, and compute the probabilities of the partitions in decreasing order of number of size 1 parts until the sum of the probabilities of the partitions that we have computed is greater than $\xi$.  The final partition computed, the partition whose probability  made the sum exceed $\xi$, is returned as $\pi$.  The {\ttfamily random} method in the {\ttfamily RandomGraph} class implements this algorithm and is available at \url{http://github.com/nicolemanno/Lights-Out}.  To create a partition of $n$ with $n-k$ parts of size 1, we create a partition of $k$ with no parts of size 1 and adjoin $n-k$ parts of size 1 to the partition.  Since the partitions of $n$ with more parts of size 1 usually have larger weight, $k$ is typically small.  In practice we are able to generate all of the partitions of $k$ in decreasing order of the largest value in the partition, $m$, by recursively partitioning the remaining $k-m$. Within the code available at our github page, the implementation of this algorithm is the {\ttfamily uniquePartitions} method within the {\ttfamily Partition} class.

\smallskip
\begin{table}[!htb]
\begin{center}
\caption{\label{tab:smallrand} For each $n$ with $1 \leq n \leq 27$, the probability that a graph chosen uniformly at random from the set of graphs with $n$ vertices is universally solvable, approximated by 1,000,000 trials. For $n \geq 11$, the margin of error at 95$\%$ confidence is less than 0.001.}
\bigskip
\begin{tabular}{ccccccccccc}
\hline\\[-5 pt]
Number & & Probability & & Probability & & Number & & Probability & & Probability \\
of & & Connected & & Universally & & of & & Connected & & Universally\\
Vertices & & & & Solvable & & Vertices & & & & Solvable\\[5 pt]
\hline\\[-5 pt]
1 && 1 && 1 && 15 && 0.999083  && 0.417295\\
2 && 0.500186 && 0.499814  && 16 && 0.999465  && 0.418368\\
3 && 0.501065 && 0.498932  && 17 && 0.999744 && 0.417804 \\
4 && 0.545260 && 0.362558 && 18 && 0.999859 && 0.419841\\
5 && 0.617609 && 0.382157  && 19 && 0.999926 && 0.419059\\
6 && 0.717691  && 0.301979 && 20 && 0.999974 && 0.419294\\
7 && 0.817183 && 0.324418  && 21 && 0.999984 && 0.420138\\
8 && 0.900342 && 0.328411 && 22 && 0.999988 && 0.419788\\
9 && 0.950290 && 0.357775 && 23 && 0.999994 && 0.419580\\
10 && 0.975863 && 0.378710 && 24 && 0.999999 && 0.419234\\
11 && 0.987933 && 0.395842  && 25 && 1 && 0.419309\\
12 && 0.993693 && 0.405955  && 26&& 0.999998 && 0.419638\\
13 && 0.996637 && 0.411057 && 27 && 0.999999 && 0.419374\\
14 && 0.998301 && 0.414513 \\[5 pt]
\hline
\end{tabular}
\end{center}
\end{table}

\subsection{Choosing Graph From Fixed Point Set}

Let $P$ be the set of two element subsets of $\{1, 2, \ldots, n\}$.  For a graph $\Gamma$ with vertices labeled from 1 to $n$, the edges $E$ of $\Gamma$ are labeled as elements of $P$.  Let $g \in \pi$ act on $P$ by $g \cdot \{i,j\} = \{g(i),g(j)\}$. For $\Gamma$ $\in \fix(g)$, each orbit in $P$ under the action of $g$ is either a subset of $E$ or disjoint from $E$. Selecting a graph $\Gamma$ uniformly at random from $\fix(g)$ reduces to computing the orbits in $P$ under the action of $g$ and constructing $\Gamma$ by giving each orbit probability 0.5 of being in $E$. We compute the orbits in $P$ under $g$ by exhaustively applying $g$ to each element of $P$. The code for this algorithm is available at our github page, and is implemented by the method {\ttfamily orbits} in the {\ttfamily RandomGraph} class.

\smallskip
\begin{table}[!htb]
\begin{center}
\caption{\label{tab:smallrandconn} For each $n$ with $1 \leq n \leq 27$, the probability that a graph chosen uniformly at random from the set of connected graphs with $n$ vertices is universally solvable, approximated by 1,000,000 trials. For $n \geq 11$, the margin of error at 95$\%$ confidence is less than 0.001.}
\bigskip
\begin{tabular}{ccccccc}
\hline\\[-5 pt]
Number & &  Probability & & Number & &  Probability\\
of & &  Universally & & of & &  Universally\\
Vertices & &  Solvable & & Vertices & &  Solvable\\[5 pt]
\hline\\[-5 pt]
1 && 1 && 15 && 0.417291 \\
2 && 0 && 16 && 0.418391 \\
3 && 0 .499886 && 17 && 0.417805 \\
4 && 0.332844 && 18 &&  0.419837  \\
5 && 0.428600 && 19 && 0.419054 \\
6 &&  0.294965 && 20 && 0.419296 \\
7 &&  0.340211 && 21 && 0.420136 \\
8 &&  0.333021 && 22 && 0.419794 \\
9 && 0.360898 && 23 && 0.419580 \\
10 && 0.379616 && 24 && 0.419235 \\
11 && 0.396218 && 25 && 0.419309\\
12 && 0.406112 && 26 && 0.419638 \\
13 && 0.411093 && 27 && 0.419373 \\
14 && 0.414537 && && \\[5 pt]
\hline
\end{tabular}
\end{center}
\end{table}

\subsection{Results}

Following the procedure given by Dixon and Wilf in \cite{Dixon}, our algorithms to select a partition $\pi$ and graph $\Gamma$, from Subsections 3.1 and 3.2, respectively, produce a representative labeled graph $\Gamma$ of an unlabeled graph chosen uniformly from the set of graphs with $n$ vertices. As discussed in Section \ref{sec:lin}, we check if $\Gamma$ is universally solvable by row reducing its neighborhood matrix. 

For each $n$ from 1 to 27, we ran 1,000,000 experiments where we chose a graph with $n$ vertices uniformly at random and determined if it was universally solvable. We also computed the probability of a graph with $n$ vertices being connected when chosen uniformly at random. These results are shown in Table \ref{tab:smallrand}.  For those same values of $n$, we also ran 1,000,000 experiments where we chose a connected graph with $n$ vertices uniformly at random and determined if it was universally solvable.  We accomplished this by choosing a graph with $n$ vertices uniformly at random and re-running the experiment if the chosen graph was not connected.  We checked if the graph is connected by starting at node 1 and performing a depth first traversal of the spanning tree of the graph made by visiting the previously unvisited adjacent vertex with smallest numerical label. If every vertex in the graph had been visited by the end of the traversal, then the graph was connected.  Within the code available at our github page, the {\ttfamily traverse} and {\ttfamily isConnected} methods in the {\ttfamily MatrixGenerator} class implement this algortihm.  Our results for connected graphs with $n$ from 1 to 27 are shown in Table \ref{tab:smallrandconn}.  For each $n$ from 28 to 100, we did not handle the connected case separately.  As $n$ approaches $\infty$, the probability that the chosen graph is connected approaches 1.  For each of these values of $n$, all 1,000,000 graphs chosen uniformly at random were connected, and we determined if each of the 1,000,000 chosen graphs was universally solvable.  We present our results for these values in Table \ref{tab:largerand}. 

For values of $n \geq 11$, our sample size is small relative to the set of all graphs with $n$ vertices; we can assume that our experiment produces a binomial distribution giving a margin of error of less than 0.001 at 95$\%$ confidence.  As confirmation of the correctness of our algorithms and implementation, observe that for values of $n$ from 1 to 11, our results in Tables \ref{tab:smallrand} and \ref{tab:smallrandconn} closely match the results in Section \ref{sec:lin}.  Additionally, for each $n$ from 1 to 27, the probability that a graph chosen uniformly at random from the set of graphs with $n$ vertices is connected, as given in Table \ref{tab:smallrand}, agrees with known values \cite{Briggs}.

\begin{table}[!htb]
\begin{center}
\caption{ \label{tab:largerand}
For each $n$ with $28 \leq n \leq 100$, the probability that a graph chosen uniformly at random from the set of graphs with $n$ vertices is universally solvable, approximated by 1,000,000 trials; all graphs chosen were incidentally connected. For each $n$, the margin of error at 95$\%$ confidence is less than 0.001.}
\bigskip
\begin{tabular}{cccccc}
\hline\\[-5 pt]
Number & Probability & Number & Probability & Number & Probability \\
of & Universally & of & Universally & of & Universally \\
Vertices &Solvable & Vertices & Solvable & Vertices & Solvable \\[5 pt]
\hline\\[-5 pt]
28 & 0.419273  & 53 & 0.419991  & 78 & 0.419690   \\
29 & 0.419329  & 54 & 0.419592 & 79 & 0.418720 \\
30 & 0.419358  & 55 & 0.420166  & 80 & 0.419922  \\
31 & 0.419409 & 56 & 0.420290  & 81 & 0.419727  \\
32 & 0.419396 & 57 & 0.418865  & 82 & 0.419339  \\
33 & 0.418807  & 58 & 0.419681  & 83 & 0.419351 \\
34 & 0.419929  & 59 & 0.419318 & 84 & 0.420456  \\
35 & 0.418619  & 60 & 0.419115 & 85 & 0.419437 \\
36 & 0.419861  & 61 & 0.418608  & 86 & 0.419870  \\
37 & 0.420124  & 62 & 0.418717  & 87 & 0.419644\\
38 & 0.419721 & 63 & 0.419397 & 88 & 0.418764 \\
39 & 0.419141  & 64 & 0.419275  & 89 & 0.419063\\
40 & 0.420303  & 65 & 0.419000 & 90 & 0.419413 \\
41 & 0.419691  & 66 & 0.419709  & 91 & 0.418826 \\
42 & 0.419005 & 67 & 0.419992  & 92 & 0.419806  \\
43 & 0.419342  & 68 & 0.418893  & 93 & 0.420023 \\
44 & 0.419507  & 69 & 0.419053  & 94 & 0.420860 \\
45 & 0.419345  & 70 & 0.420118  & 95 & 0.420415\\
46 & 0.420292 & 71 & 0.419159 & 96 & 0.418757  \\
47 & 0.419140 & 72 & 0.419537  & 97 & 0.419563 \\
48 & 0.418924  & 73 & 0.418961  & 98 & 0.419010 \\
49 & 0.418542  & 74 & 0.419708 & 99 &0.419510  \\
50 & 0.418690 & 75 & 0.419446  & 100 & 0.419362\\
51 & 0.419030 & 76 & 0.419571\\
52 & 0.419023 & 77 & 0.419299\\[5 pt]
\hline
\end{tabular}
\end{center}
\end{table}

The runtime for the results displayed in Tables \ref{tab:smallrand} and \ref{tab:smallrandconn} was 30 minutes, while it took 48 hours to compute the results for Table \ref{tab:largerand}.  These computations were performed on a laptop with a 2.80 GHz processor and 16 GB of RAM running Windows 11.

\section{Future Work}
\label{sec:future}
There are a number of potential extensions of this work. Variants of \emph{Lights Out} with more than one toggle mode have been considered by many authors \cite{BBS} \cite{EEJJMS} \cite{GP} \cite{GMT} \cite{HMP} \cite{Zaidenberg}. It is possible to extend our results in this direction.

Instead of choosing a graph uniformly at random from the set graphs with $n$ vertices, one could restrict the set of graphs under consideration to those with $m$ edges and $n$ vertices.  Fixing $n$ and examining how the probability of choosing a universally solvable graph changes as $m$ changes is an interesting direction to extend this work.

Our work suggests that as $n$ approaches $\infty$, the probability that a graph with $n$ vertices chosen uniformly at random is universally solvable is approximately $0.419$. This is only an approximation and conjecture, investigating this is an intriguing open problem.

\section*{Acknowledgments}

The authors acknowledge the contributions of Noah Williams and the anonymous referee.  Their insightful suggestions and helpful comments improved this work significantly.

{\footnotesize


\begin{thebibliography}{00}

\bibitem{AS1} A.T. Amin, P.J. Slater, Neighborhood domination with parity restrictions in graphs, {\it Congr. Numer.}, {\bf 91} (1992), 19--30. \MR{1208985}

\bibitem{AS2} A.T. Amin, P.J. Slater, All parity realizable trees, {\it J. Combin. Math. Combin. Comput.}, { \bf 20} (1996), 53--63. \MR{1376697}

\bibitem{Anderson} M. Anderson, T. Feil, \href{https://doi.org/10.1080/0025570X.1998.11996658}{Turning lights out with linear algebra}, { \it Math. Mag.}, { \bf 71} (1998), 300--303. \MR{1573341}

\bibitem{BBS} L. Ballard, E. Budge, D. Stephenson, \href{https://doi.org/10.2140/involve.2019.12.181}{Lights out for graphs related to one another by constructions},  { \it Involve}, { \bf 12} (2019), 181--201. \MR{3864213}

\smallskip
\bibitem{Briggs}
K. Briggs, Combinatorial Graph Theory, (2022, May 11). Retrieved from \url{http://keithbriggs.info/cgt.html}

\bibitem{CHKR} R. Cowen, S.H. Hechler, J.W. Kennedy, A. Ryba, Inversion and neighborhood inversion in graphs, { \it Graph Theory Notes N.Y.}, { \bf 37} (1999), 37--41. \MR{1725188}

\bibitem{Dixon} J.D. Dixon, H.S. Wilf, \href{https://doi.org/10.1016/0196-6774(83)90021-4}{The random selection of unlabeled graphs}, { \it J. Algorithms}, { \bf 4} (1983), 205--213. \MR{0710717}

\bibitem{EEJJMS} S. Edwards, V. Elandt, N. James, K. Johnson, Z. Mitchell, D. Stephenson, \href{https://doi.org/10.2140/involve.2010.3.17}{Lights Out on finite graphs}, { \it Involve}, { \bf 3} (2010), 17--32. \MR{2672500}

\bibitem{EES} H. Eriksson, K. Eriksson, J. Sj\"ostrand, \href{https://doi.org/10.1006/aama.2001.0739}{Note on the lamp lighting problem}, { \it Adv. Appl. Math.}, { \bf 27} (2004), 357--366 \MR{1868970}

\bibitem{GP} A. Giffen, D.B. Parker, On generalizing the \emph{Lights Out} game and a generalization of parity domination, { \it Ars Combin.}, { \bf 111} (2013), 273--288. \MR{3100179}

\bibitem{GKT} J. Goldwasser,  W. Klostermeyer, G. Trapp, \href{https://doi.org/10.1080/03081089708818520}{Characterizing switch-setting problems}, { \it Linear Multilinear Algebra}, { \bf 43} (1997), 121--136.  \MR{1613183}

\bibitem{GKTZ} J. Goldwasser, W. Klostermeyer, G. Trapp, C. Zhang, {\it Setting switches in a grid, Technical Report TR-95-20}, Department of Statistics and Computer Science, West Virginia University, 1995.

\bibitem{GMR} C. Gray, L. Mitchell, M. Roughan, \href{https://doi.org/10.1093/comnet/cnz011}{Generating connected random graphs}, { \it J. Complex Netw.}, { \bf 7} (2019), 896--912. \MR{4041485}

\bibitem{GMT} S. Gravier, M. Mhalla, E. Tannier, \href{https://doi.org/10.1016/S0304-3975(03)00285-8}{On a modular domination game}, { \it Theoret. Comput. Sci.}, { \bf 306} (2003), 291--303. \MR{2000178} 

\bibitem{HY} Y. Hayata, M. Yamagishi,\href{https://doi.org/10.2140/involve.2019.12.713}{ On weight-one solvable configurations of Lights Out Puzzle}, { \it Involve}, { \bf 12} (2019), 713--720. \MR{3941607}  

\bibitem{HMP} M. Hunziker, A. Machiavelo, J. Park,\href{https://doi.org/10.1016/j.tcs.2004.03.031}{ Chebyshev polynomials over finite fields and reversibility of $\sigma$-automata on square grids}, { \it Theoret. Comput. Sci.}, { \bf 320} (2004), 465--483. \MR{2064312}

\bibitem{McKay1} B.D. McKay, Combinatorial Data, available online at the URL: \url{http://users.cecs.anu.edu.au/~bdm/data/graphs.html}

\bibitem{McKay2} B.D. McKay, A. Piperno, \href{http://dx.doi.org/10.1016/j.jsc.2013.09.003}{Practical Graph Isomorphism, II}, { \it J. Symbolic Comput.}, { \bf 60} (2013), 94--112. \MR{3131381} 

\bibitem{Ober} W. Oberschelp, \href{https://doi.org/10.1007/BF01363123}{Kombinatorische Anzahlbestimmungen in Relationen}, { \it Math. Ann.}, { \bf 174} (1967), 53--78.  \MR{0218255}

\bibitem{oeis} N.J. Sloane, A000088, available online at the URL: \url{http://oeis.org/A000088}

\bibitem{Sutner1} K. Sutner, \href{https://doi.org/10.1007/BF03023823}{Linear cellular automata and the Garden-of-Eden}, { \it Math. Intelligencer}, { \bf 11} (1989), 49--53. \MR{0994964} 

\bibitem{Sutner2} K. Sutner, \href{https://doi.org/10.2307/2323999}{The $\sigma$-game and cellular automata}, { \it Amer. Math. Monthly}, { \bf 97} (1990), 24--34. \MR{1034347}  

\bibitem{Zaidenberg} M. Zaidenberg, \href{https://doi.org/10.2478/s11533-009-0029-0}{Periodic harmonic functions on lattices and points count in positive characteristic}, { \it Cent. Eur. J. Math.}, { \bf 7} (2009), 365--381. \MR{2534458}

\end{thebibliography}
\end{document}